# CHARACTERIZATION OF ARBITRAGE-FREE MARKETS[1]

By Eva Strasser

*Vienna University of Technology*

The present paper deals with the characterization of no-arbitrage properties of a continuous semimartingale. The first main result, Theorem 2.1, extends the no-arbitrage criterion by Levental and Skorohod [*Ann. Appl. Probab.* **5** (1995) 906–925] from diffusion processes to arbitrary continuous semimartingales. The second main result, Theorem 2.4, is a characterization of a weaker notion of no-arbitrage in terms of the existence of supermartingale densities. The pertaining weaker notion of no-arbitrage is equivalent to the absence of immediate arbitrage opportunities, a concept introduced by Delbaen and Schachermayer [*Ann. Appl. Probab.* **5** (1995) 926–945].

Both results are stated in terms of conditions for any semimartingales starting at arbitrary stopping times $\sigma$. The necessity parts of both results are known for the stopping time $\sigma = 0$ from Delbaen and Schachermayer [*Ann. Appl. Probab.* **5** (1995) 926–945]. The contribution of the present paper is the proofs of the corresponding sufficiency parts.

**1. Introduction.** In a discrete-time model, the usual definition of the no-arbitrage property (NA-property) can be characterized by the existence of an equivalent martingale measure for the underlying price process; see Harrison and Pliska (1981) and Dalang, Morton and Willinger (1990). Within the setting of a continuous-time model, Kreps (1981) associates the existence of an equivalent martingale measure with a stronger no-arbitrage property, the so-called property of no free lunch (NFL-property). In a series of detailed studies, Delbaen and Schachermayer (1994, 1995) show that the NFL-property is equivalent to the apparently weaker property of no free lunch with vanishing risk (NFLVR-property), clarifying the situation with various versions of the NA-property.

Received September 2002; revised January 2004.

[1]Supported by the Austrian Science Foundation (FWF) under the Wittgenstein–Preis program Z36-MAT and Grant Spezialforschungsbereich SFB#0101.

The focus of this paper is on the characterization of weaker notions of no-arbitrage than the NFL-property and the NFLVR-property, respectively. In this regard, we extend criteria going back to Levental and Skorohod (1995) and Schweizer (1995), applying results by Delbaen and Schachermayer (1995).

Fix a finite time horizon $T > 0$ and a stochastic base $(\Omega, \mathcal{F}_T, P; \mathbf{F})$, where the filtration $\mathbf{F} = (\mathcal{F}_t)_{0 \leq t \leq T}$ is assumed to satisfy the usual conditions. The set of $\mathbb{R}^d$-valued semimartingales is denoted by $\mathcal{S}$ and for $S \in \mathcal{S}$, the set $L(S)$ is defined to be the set of $\mathbb{R}^d$-valued, predictable, $S$-integrable processes. Moreover, recall that a locally square integrable semimartingale has a canonical decomposition, $S = S_0 + M + A$, into a locally square integrable martingale $M$ with $M_0 = 0$ and a predictable process of finite variation $A$ with $A_0 = 0$.

In the context of mathematical finance, a market model is a vector of $d+1$ assets, one bond and $d$ stocks. The price process of the bond is assumed to be constant (i.e., we choose the bond as numeraire) and the price process of the $d$ stocks is assumed to be an $\mathbb{R}^d$-valued semimartingale $S$. A portfolio is a pair $(x, H)$, where $x \in \mathbb{R}$ is the initial wealth and $H \in L(S)$ specifies the number of shares in each asset held in the portfolio. The corresponding (self-financing) wealth-process $X$ is given by $X = x + H \cdot S$. Let us denote by $\mathcal{X}$ the family of wealth-processes, that is,

$$\mathcal{X} := \{X = x + H \cdot S : x \in \mathbb{R}, H \in L(S)\}.$$

A wealth-process $X \in \mathcal{X}$ is called admissible if it is uniformly bounded from below. Let us denote by $\mathcal{X}_a \subseteq \mathcal{X}$ the family of admissible wealth-processes and by $\mathcal{X}_+ \subseteq \mathcal{X}_a$ the family of nonnegative wealth-processes.

For a semimartingale $S$ and a stopping time $\sigma \leq T$ we define $^\sigma S$ to be the semimartingale $S$ starting at $\sigma$, that is, $^\sigma S_t := S_{(\sigma+t) \wedge T}$, $0 \leq t \leq T$. Note that $^\sigma S$ is adapted to the filtration $^\sigma \mathbf{F} = (^\sigma \mathcal{F}_t)_{0 \leq t \leq T}$, where $^\sigma \mathcal{F}_t := \mathcal{F}_{(\sigma+t) \wedge T}$, $0 \leq t \leq T$. Accordingly, we define $^\sigma \mathcal{X}$ to be the set of wealth-processes given by $^\sigma X = x + H \cdot {}^\sigma S$, where $x \in \mathbb{R}$ and $H \in L(^\sigma S)$.

Let us recall some basic concepts of no-arbitrage theory. We say that $S$ satisfies the NA-property, if for every $X \in \mathcal{X}_a$, we have

(NA)　　　　　　　$X_0 = 0$　and　$X_T \geq 0$　$\implies$　$X_T = 0$.

Moreover, we say that $S$ satisfies the $\mathrm{NA}^+$-property, if for every $X \in \mathcal{X}_+$, we have

($\mathrm{NA}^+$)　　　　　　$X_0 = 0$　and　$X \geq 0$　$\implies$　$X = 0$.

Recall that by Delbaen and Schachermayer [(1995), Lemma 3.1], for a continuous semimartingale the $\mathrm{NA}^+$-property is equivalent to the absence of so-called immediate arbitrage opportunities.

An equivalent (absolutely continuous) probability measure $Q \sim P$ ($Q \ll P$) is called an equivalent (absolutely continuous) local martingale measure for the semimartingale $S$ if $S$ is a local $Q$-martingale with respect to the filtration $\mathbf{F}$.

**2. Main results.** It is a well-known fact that the existence of an equivalent local martingale measure implies the NA-property. The reverse implication is not true in general. However, for special classes of semimartingales one can achieve a characterization. This is the topic of our first main result, Theorem 2.1.

Levental and Skorohod (1995) prove in their Theorem 1 a kind of prototype of our Theorem 2.1 under the additional assumption that the local martingale part $M$ of the continuous semimartingale $S$ is of the form $M = \Sigma \cdot W$. Here, $W$ is an $\mathbb{R}^d$-valued Brownian motion defined on its natural filtration and $\Sigma$ is an adapted matrix-valued process such that each $\Sigma_t$ is invertible, $0 \leq t \leq T$. Within this framework, the martingale representation property holds true and thus the proof can be based on explicit representations of the local martingale measures. Delbaen and Schachermayer (1995) consider in their Theorem 1.4 the more general case of arbitrary continuous semimartingales and show that the NA-property implies the existence of an absolutely continuous local martingale measure. The proof relies on the fundamental theorem of asset pricing by Delbaen and Schachermayer (1994). Recently, Kabanov and Stricker (2003) extend this result, dropping the continuity assumption for the semimartingales.

Theorem 2.1 is an extension of the criterion by Levental and Skorohod (1995) to arbitrary continuous semimartingales, using the result by Delbaen and Schachermayer (1995). In the meantime, after the submission of the present paper, Kabanov and Stricker (2003) extended Theorem 2.1 to the case of markets with countably many assets.

THEOREM 2.1. *The continuous semimartingale $S$ satisfies the NA-property iff for every stopping time $\sigma \leq T$ there exists an absolutely continuous local martingale measure $^\sigma Q \ll P$ satisfying $^\sigma Q|_{\mathcal{F}_\sigma} \sim P|_{\mathcal{F}_\sigma}$ for the semimartingale $^\sigma S$ starting at $\sigma$.*

For a detailed proof we refer to Section 3. The author thanks Y. Kabanov for pointing out a lacuna in a preliminary version of the proof.

Before we present our second main result, Theorem 2.5, we give a reformulation of Theorem 1 in Levental and Skorohod (1995). For this purpose, let us recall an important structure condition, which can be characterized in terms of a very weak notion of no-arbitrage.

THEOREM 2.2. *Let $S$ be a locally square integrable semimartingale. Then every nonnegative, predictable wealth-process $X \in \mathcal{X}$ of bounded variation is constant iff the structure condition $dA \ll d\langle M, M \rangle$ is valid, that is, there exists a predictable process $\lambda$ with values in $\mathbb{R}^d$ such that $dA = d\langle M, M\rangle \lambda$.*

This theorem is proved using the same arguments as in the proof of Theorem 3.5 in Delbaen and Schachermayer (1995). A direct proof is given in Strasser (2003).

A semimartingale $S$ as used by Levental and Skorohod (1995) always satisfies the structure condition $dA \ll d\langle M, M\rangle$. In this case, the absolutely continuous local martingale measures used in Theorem 2.1 can be derived from a particular process. To be explicit, such a semimartingale $S$ satisfies the NA-property iff for every stopping time $\sigma \leq T$ the density process

$$(2.1) \qquad {}^{\sigma}Z := \mathcal{E}\left(-\int \lambda' I_{]\!] \sigma, T]\!]}\, dM\right)$$

satisfies $E({}^{\sigma}Z_T) = 1$. The relation between our Theorem 2.1 and Levental and Skorohod (1995), Theorem 1, is then established by defining $d\mathcal{Q} := {}^{\sigma}Z_T dP$.

Now we are in a position to present our second main result. Recall that for an arbitrary continuous semimartingale $S$ satisfying $dA \ll d\langle M, M\rangle$, the condition $E({}^{\sigma}Z_T) = 1$ implies the existence of an absolutely continuous local martingale measure for ${}^{\sigma}S$. In general, if $E({}^{\sigma}Z_T) < 1$, ${}^{\sigma}Z$ is not the density of an absolutely continuous local martingale measure, but sometimes a so-called supermartingale density, a notion we adopt from Kramkov and Schachermayer [(1999), Section 2]. Note that this notion of a supermartingale density is slightly stronger than that introduced by Kabanov and Stricker (2003).

DEFINITION 2.3. A nonnegative (strictly positive) process $Y$ is a (strict) supermartingale density for $S$ if $Y$ is a supermartingale with $Y_0 = 1$ such that the product $YX$ is a supermartingale for every $X \in \mathcal{X}_+$.

It is easy to see that the existence of a strict supermartingale density implies the $NA^+$-property. Our second main theorem shows that it is even possible to characterize the $NA^+$-property by means of the existence of supermartingale densities. The key idea is similar to that of Theorem 2.1, that is, we impose conditions on the semimartingales ${}^{\sigma}S$ starting at arbitrary stopping times $\sigma$. For a proof as well as further equivalent assertions see Theorem 3.5 in Section 3.

THEOREM 2.4. *A continuous semimartingale $S$ satisfies the $NA^+$-property iff for every stopping time $\sigma \leq T$ the process ${}^{\sigma}Z$ defined in (2.1) is a supermartingale density for ${}^{\sigma}S$.*

**3. Proof of the main results.** Let us begin with the proof of Theorem 2.1. For this purpose, we state and prove an auxiliary lemma, which isolates the basic idea of the proof of Theorem 2.1.

LEMMA 3.1. *Let $S$ be a continuous semimartingale and suppose that there exists an absolutely continuous local martingale measure $Q \ll P$ satisfying $Q|_{\mathcal{F}_0} \sim P|_{\mathcal{F}_0}$. Let $X \in \mathcal{X}_a$ such that $X_0 = 0$ and $X_T \geq 0$ and define*

$$\tau := \inf\{t > 0 : X_t \neq 0\}.$$

*Then we have $\tau > 0$.*

PROOF. Denote by $Z$ the density process of $Q$ with respect to $P$. The local martingale $ZX$ is bounded from below by a multiple of the martingale $Z$ and thus $ZX$ is a supermartingale. The nonnegativity of $Z_T X_T$ yields $ZX \equiv 0$ and thus $X^\theta \equiv 0$, where the stopping time $\theta$ is defined by

$$\theta := \inf\{t > 0 : Z_t = 0 \text{ or } Z_{t-} = 0\}$$

satisfying $\theta > 0$. For a more detailed proof we refer to Strasser (2003). □

PROOF OF THEOREM 2.1. *Necessity*: Fix a stopping time $\sigma \leq T$ and observe that the NA-property of $S$ implies the NA-property of $^\sigma S$. Applying Delbaen and Schachermayer [(1995), Theorem 1.4], we get the existence of an absolutely continuous local martingale measure $^\sigma Q \ll P$ satisfying $^\sigma Q|_{\mathcal{F}_\sigma} \sim P|_{\mathcal{F}_\sigma}$ for $^\sigma S$. A recent discussion of this result can be found in Kabanov and Stricker (2003).

*Sufficiency*: Let $X \in \mathcal{X}_a$ with $X_0 = 0$ and $X_T \geq 0$. Define the stopping time

$$\sigma := \inf\{t > 0 : X_t \neq 0\} \wedge T$$

and assume $P(\sigma < T) > 0$. It is easy to see that $^\sigma X \in {}^\sigma \mathcal{X}_a$ satisfies $^\sigma X_0 = 0$ and $^\sigma X_T \geq 0$. Moreover,

$$^\sigma \tau := \inf\{t \geq 0 : {}^\sigma X_t \neq 0\}$$

is a stopping time with respect to the filtration $^\sigma \mathbf{F}$ satisfying $^\sigma \tau = 0$ on $\{\sigma < T\}$.

By assumption, there exists an absolutely continuous probability measure $^\sigma Q \ll P$ satisfying $^\sigma Q|_{\mathcal{F}_\sigma} \sim P|_{\mathcal{F}_\sigma}$ for $^\sigma S$. Applying Lemma 3.1 to $^\sigma S$ and $^\sigma X$ yields $^\sigma \tau > 0$ on $\{\sigma < T\}$. This is a contradiction. Hence, $\sigma = T$ and the assertion is proved. □

Let us turn to the discussion of Theorem 2.4. For this purpose, let $F$ be a càdlàg, predictable and increasing process with $F_0 = 0$ such that $dA = g\,dF$ and $d\langle M, M\rangle = v\,dF$. Here, $g$ and $v$ denote predictable Radon–Nikodym derivatives. Using this notation, the structure condition $dA \ll d\langle M, M\rangle$ can equivalently be stated as follows: there exists a predictable process $\lambda$ with values in $\mathbb{R}^d$ such that

$$g = v\lambda, \qquad F \otimes P\text{-a.e.} \tag{3.1}$$

This structure condition is frequently used in the literature; see, for example, Ansel and Stricker (1992) and Schweizer (1995). Karatzas and Shreve (1998) consider positive continuous semimartingales and naturally use a logarithmic version of (3.1).

Let us define the notion of the mean-variance trade-off, similar to Schweizer [(1995), Section 2].

DEFINITION 3.2. Assume $dA \ll d\langle M, M\rangle$ and let $\lambda$ be any predictable process satisfying condition 2 of Theorem 2.2 [or equivalently equation (3.1)]. The family

$$K_s^t := \int_s^t \lambda' v \lambda \, dF, \qquad 0 \leq s \leq t \leq T, \tag{3.2}$$

is called the mean-variance trade-off (MVT).

Clearly, the MVT is not necessarily finite. Finiteness of the MVT, that is, $K_0^T < \infty$, simplifies the situation considerably. The following assertion follows from Theorem 2.2 and Corollary 3 in Schweizer (1995).

COROLLARY 3.3. *Let $S$ be a continuous semimartingale and assume $K_0^T < \infty$. Then $S$ satisfies the $NA^+$-property iff $dA \ll d\langle M, M\rangle$.*

The question arises whether it is possible to characterize the $NA^+$-property without assuming $K_0^T < \infty$. This problem is settled by our second main result, Theorem 2.4. For completion, we will prove a more detailed assertion in Theorem 3.5 containing that of Theorem 2.4. For this we need the following notion, going back to Delbaen and Schachermayer (1995) and Levental and Skorohod (1995).

DEFINITION 3.4. We say that the MVT does not jump to $\infty$, if the stopping time $\alpha$ defined by

$$\alpha := \inf\{t > 0 : K_t^{t+h} = \infty \ \forall h \in \,]0, T-t]\}$$

satisfies $\alpha = \infty$.

Obviously, $K_0^T < \infty$ implies $\alpha = \infty$, whereas the converse is not true in general. For reasons of proof, we equivalently reformulate the assertion of Theorem 2.4 in the following theorem.

THEOREM 3.5. *Let $S$ be a continuous semimartingale. The following are equivalent*:

1. *The semimartingale $S$ satisfies the $NA^+$-property.*
2. *The structure condition $dA \ll d\langle M, M\rangle$ is valid and the MVT does not jump to $\infty$.*
3. *For every stopping time $\sigma \leq T$ the process ${}^\sigma Z$ defined in (2.1) is a supermartingale density for ${}^\sigma S$.*
4. *For every stopping time $\sigma \leq T$ there exists a supermartingale density ${}^\sigma Y$ for ${}^\sigma S$.*

The proof of the implication $1 \Rightarrow 2$ is a straightforward extension of Delbaen and Schachermayer [(1995), Section 3]. Below, we prove implications $2 \Rightarrow 3$ and $4 \Rightarrow 1$.

In the setting of Levental and Skorohod (1995), the $NA^+$-property implies $\alpha = \infty$. Our Theorem 2.4 shows that under the assumption $dA \ll d\langle M, M\rangle$, which is weaker than the setting of Levental and Skorohod (1995), the condition $\alpha = \infty$ is even equivalent to the $NA^+$-property.

PROOF OF THEOREM 3.5. $2 \Rightarrow 3$: Fix a stopping time $\sigma$ with $P(\sigma < T) > 0$. It is easy to see that $P(K_\sigma^{\sigma+h} = \infty \; \forall h \in \rrbracket 0, T-\sigma \rrbracket) = 0$ and that the stopping time

$$(3.3) \qquad \tau := \tau(c) := \inf\{h > \sigma : K_\sigma^{\sigma+h} = c\} \wedge T, \qquad c \in \mathbb{R}_+,$$

satisfies $\tau > \sigma$ on $\{\sigma < T\}$, since the MVT $(K_\sigma^{\sigma+h})_{h \in \rrbracket 0, T-\sigma \rrbracket}$ is continuous, starts at zero and does not jump to $\infty$. Moreover, $\int \lambda' v \lambda \mathbb{1}_{\rrbracket \sigma, \tau \rrbracket} \, dF < \infty$ and thus $\int \lambda' \mathbb{1}_{\rrbracket \sigma, \tau \rrbracket} \, dM$ is a locally square integrable martingale.

Note that the continuous semimartingale ${}^\sigma S$ has a canonical decomposition ${}^\sigma S = {}^\sigma S_0 + {}^\sigma M + {}^\sigma A$ with respect to the filtration ${}^\sigma \mathbf{F}$. Moreover, the predictable characteristics ${}^\sigma g$ and ${}^\sigma v$ of ${}^\sigma S$ satisfy ${}^\sigma g = (g_{\sigma+t}^T)_{t \in [0,T]}$ and ${}^\sigma v = (v_{\sigma+t}^T)_{t \in [0,T]}$. Finally, $dA \ll d\langle M, M\rangle$ obviously implies $d{}^\sigma A \ll d\langle {}^\sigma M, {}^\sigma M\rangle$ and thus we have ${}^\sigma g = {}^\sigma v {}^\sigma \lambda$, $P \otimes F$-a.e.

Define $\rho := \rho(c) := \tau - \sigma$, where $\tau$ is as in (3.3), and note that $\rho \geq 0$ as well as $\rho > 0$ on $\{\sigma < T\}$. In general, $\rho$ is not a stopping time with respect to the filtration $\mathbf{F}$, but it is a stopping time with respect to the filtration ${}^\sigma \mathbf{F}$. Indeed, since $\sigma$ and $\tau$ are stopping times with respect to the filtration $\mathbf{F}$, we obtain

$$\{\rho \leq t\} = \{\tau - \sigma \leq t\}$$
$$= \{\tau \leq (\sigma + t) \wedge T\} \in \mathcal{F}_\tau \cap \mathcal{F}_{(\sigma+t) \wedge T} \subseteq {}^\sigma \mathcal{F}_t \qquad \forall t \in [0, T].$$

Define

$$(3.4) \qquad {}^\sigma Z := \mathcal{E}\left(-\int {}^\sigma \lambda' \, d{}^\sigma M\right).$$

Since $\int {}^\sigma\lambda' \, d({}^\sigma M)^\rho$ is a locally square integrable martingale, it follows that $0 < ({}^\sigma Z)^\rho = 1 - \int {}^\sigma Z {}^\sigma \lambda' \, d({}^\sigma M)^\rho$ is a locally square integrable martingale, too, and that

$$(3.5) \quad ({}^\sigma A)^\rho = \int {}^\sigma g I_{\rrbracket 0,\rho \rrbracket} \, d^\sigma F = \int {}^\sigma v {}^\sigma \lambda I_{\rrbracket 0,\rho \rrbracket} \, d^\sigma F = \int {}^\sigma \lambda' \, d\langle {}^\sigma M, {}^\sigma M \rangle^\rho.$$

Straightforward computations as in Schweizer (1995) prove the local martingale property of $({}^\sigma Z {}^\sigma S)^\rho$. Consequently, $({}^\sigma Z {}^\sigma X)^\rho$ is a supermartingale for every ${}^\sigma X \in {}^\sigma \mathcal{X}_+$.

Define $\rho(\infty) := \tau(\infty) - \sigma$, where

$$\tau(\infty) := \inf\{h > \sigma : K_\sigma^{\sigma+h} = \infty\} \wedge T.$$

Observe $\lim_{c \to \infty} \rho(c) = \rho(\infty)$ and ${}^\sigma Z \mathbb{1}_{\rrbracket \rho(\infty), T \rrbracket} = 0$. In particular, we have ${}^\sigma Z = ({}^\sigma Z)^{\rho(\infty)}$ as well as ${}^\sigma Z {}^\sigma X = ({}^\sigma Z {}^\sigma X)^{\rho(\infty)}$. Together with Fatou's lemma, this implies

$$E({}^\sigma Z_t {}^\sigma X_t | \mathcal{F}_s) \leq \lim_{c \to \infty} E(({}^\sigma Z_t {}^\sigma X_t)^{\rho(c)} | \mathcal{F}_s) \leq \lim_{c \to \infty} ({}^\sigma Z_s {}^\sigma X_s)^{\rho(c)} = {}^\sigma Z_s {}^\sigma X_s,$$

$0 \leq s \leq t \leq T$, since the choice of $c \in \mathbb{R}_+$ in (3.3) was arbitrary. Hence, ${}^\sigma Z$ is a supermartingale density for ${}^\sigma S$.

$4 \Rightarrow 1$: Let $X \in \mathcal{X}_+$ with $X_0 = 0$. Define the stopping time

$$(3.6) \qquad \sigma := \inf\{t > 0 : X_t \neq 0\} \wedge T$$

and assume $P(\sigma < T) > 0$. It is easy to see that ${}^\sigma X \in {}^\sigma \mathcal{X}_+$ satisfies ${}^\sigma X_0 = 0$. Moreover,

$${}^\sigma \tau := \inf\{t \geq 0 : {}^\sigma X_t \neq 0\} \wedge T$$

is a stopping time with respect to the filtration ${}^\sigma \mathbf{F}$ satisfying ${}^\sigma \tau = 0$ on $\{\sigma < T\}$.

By assumption, there exists a supermartingale density ${}^\sigma Y$ for ${}^\sigma S$. Define the stopping time

$$\theta := \inf\{t \geq 0 : {}^\sigma Y_t = 0 \text{ or } {}^\sigma Y_{t-} = 0\}$$

and note that $\theta > 0$ on $\{\sigma < T\}$ since ${}^\sigma Y_0 = 1$. Since ${}^\sigma Y {}^\sigma X$ is a nonnegative supermartingale, it follows that ${}^\sigma Y {}^\sigma X \equiv 0$ and $({}^\sigma X)^\theta \equiv 0$. We obtain in particular that ${}^\sigma \tau \geq \theta > 0$ on $\{\sigma < T\}$. This is a contradiction. Hence, $\sigma = T$ and the assertion is proved. $\square$

**Acknowledgments.** The author thanks F. Delbaen, Y. Kabanov and W. Schachermayer for discussions on the topic of the present paper as well as valuable suggestions. Special thanks go to the referee whose suggestions improved the presentation considerably.

# REFERENCES

Ansel, J. P. and Stricker, C. (1992). Lois de martingale, densités et décomposition de Föllmer–Schweizer. *Ann. Inst. H. Poincaré Probab. Statist.* **28** 375–392. MR1183992

Dalang, R. C., Morton, A. and Willinger, W. (1990). Equivalent martingale measures and no-arbitrage in stochastic securities market model. *Stochastics Stochastics Rep.* **29** 185–201. MR1041035

Delbaen, F. and Schachermayer, W. (1994). A general version of the fundamental theorem of asset pricing. *Math. Ann.* **300** 463–520. MR1304434

Delbaen, F. and Schachermayer, W. (1995). The existence of absolutely continuous local martingale measures. *Ann. Appl. Probab.* **5** 926–945. MR1384360

Harrison, M. and Pliska, S. (1981). Martingales and stochastic integrals in the theory of continuous trading. *Stochastic Process. Appl.* **11** 215–260. MR622165

Kabanov, Y. and Stricker, Ch. (2003). Remarks on the true no-arbitrage property. *Laboratoire de Mathématiques de Besançon.*

Karatzas, I. and Shreve, S. (1998). *Methods of Mathematical Finance.* Springer, New York. MR1640352

Kramkov, D. and Schachermayer, W. (1999). The asymptotic elasticity of utility functions and optimal investment in incomplete markets. *Ann. Appl. Probab.* **9** 904–950. MR1722287

Kreps, D. M. (1981). Arbitrage and equilibrium in economies with infinitely many commodities. *J. Math. Econom.* **8** 15–350. MR611252

Leventpal, S. and Skorohod, A. V. (1995). A necessary and sufficient condition for absence of arbitrage with tame portfolios. *Ann. Appl. Probab.* **5** 906–925. MR1384359

Schweizer, M. (1995). On the minimal martingale measure and the Föllmer–Schweizer decomposition. *Stochastic Anal. Appl.* **13** 573–599. MR1353193

Strasser, E. (2003). No-arbitrage and utility maximization in mathematical finance. Ph.D. thesis.

Department of Financial
 and Actuarial Mathematics
Vienna University of Technology
Wiedner Hauptstrasse 8-10/105
1040 Vienna
Austria
e-mail: evastrasser@fam.tuwien.ac.at
url: www.fam.tuwien.ac.at/evastrasser